\RequirePackage{fix-cm}
\documentclass{svjour3}                     
\usepackage{graphicx}
\usepackage{latexsym,amssymb,amsfonts,amsmath}
\usepackage{mathrsfs}
\usepackage[usenames]{color}
\usepackage{eucal}

\usepackage{tikz}
\usetikzlibrary{decorations.markings}
\usepackage{pgfplots}

\newcommand{\Real}{\mathbb R}
\newcommand{\eps}{\varepsilon}

\newcommand{\F}{\mathcal{F}}

\newcommand{\one}[1]{\mathbf{1}_{\{#1\}}}
\newcommand{\oneset}[1]{\mathbf{1}_{#1}}
\renewcommand{\P}{\mathbb{P}}
\newcommand{\E}{\mathbb{E}}

\journalname{J. Math. Biol.}
\pgfplotsset{compat=1.12}
\begin{document}

\title{What can be observed in real time PCR and when does it show?}

%


\titlerunning{How did the PCR start?}        

\author{Pavel Chigansky \and Peter Jagers \and Fima C.\ Klebaner
\thanks{Research  supported   by Australian Research Council Grant
No  DP150103588. P. Jagers was further supported by the Knut och Alice Wallenberg Foundation.
P. Chigansky was supported by ISF 558/13 grant}}


\institute{P. Chigansky \at
Department of Statistics,
The Hebrew University,
Mount Scopus, Jerusalem 91905,
Israel
\email{pchiga@mscc.huji.ac.il}
 \and
 P. Jagers \at Mathematical Sciences, Chalmers University of Technology
and University of Gothenburg, SE-412 96 Gothenburg, Sweden
\email{jagers@chalmers.se}
\and 
F. C. Klebaner
\at School of Mathematical Sciences,
Monash University, Monash, VIC 3800, Australia
\email{fima.klebaner@monash.edu}
}

\date{Received: date / Accepted: date}

\maketitle

\begin{abstract}
Real time, or quantitative, PCR typically starts from a very low concentration of initial DNA strands. During iterations the numbers increase, first essentially by doubling, later predominantly in a linear way. Observation of the number of DNA molecules in the experiment becomes possible only when
it is substantially larger than initial numbers, and then possibly affected by the randomness in individual replication. Can the initial copy number still be determined? 

This is a classical problem and, indeed, a concrete special case of the general problem of determining the number of ancestors, mutants or invaders, of a population observed only later. 
We approach it through a generalised version of the
branching process model introduced in \cite{JagKlePCR}, and based on Michaelis-Menten type enzyme kinetical considerations from \cite{SchnellM}. A crucial role is played by the Michaelis-Menten constant being large, as compared to initial
copy numbers. In a strange way, determination of the initial number turns out to be completely possible if the initial rate $v$ is one, i.e all DNA strands replicate, but only partly so when $v<1$, and thus the initial rate or probability of succesful replication is lower than one. Then, the starting molecule number becomes hidden behind a ``veil of uncertainty''.  This is a special  case, of a hitherto unobserved general phenomenon in population growth processes, which will be adressed elsewhere.


\end{abstract}

\nocite{*}

\section{Introduction}
In the polymerase chain reaction a molecule replicates with a
probability $p(z)$, which will be of the form
\[p(z)=\frac{C}{K+z},\]
under the asumption of Michaelis-Menten kinetics. Here, $K$ is the
Michaelis-Menten constant, large in terms of molecule numbers, $z$ the number of DNA molecules at the
actual round, and $C$ a constant, which can be written as $vK$, 
where $v$ is the maximal rate or speed of the reaction, corresponding to
$z=0$. Then, $v = p(0)$ is the probability of successful replication
under the most benign circumstances, and the decrease of $p(z)$, as
the number $z$ of DNA strands present increases, mirrors that the
latter are being synthesized from DNA building blocks, which disappear
as the number of DNA molecules increases. As has been observed
recently, though this is the general pattern, there are exceptions
where the replication probability actually increases in the very first
generation,  due to impurities in templates \cite{oralst}.

In this paper we disregard this and rely upon the Michaelis-Menten
based approach in \cite{JagKlePCR}, where it was used to
explain the first exponential but later linear growth of molecule
numbers,  see also \cite{Best} \cite{Lalam04}, \cite{Lievens}. 
For a statistical analysis, where PCR is modeled by branching
processes without environmental change due to 
growth but with random effects and starting numbers cf. \cite{Vidyashankar}.

Here we turn to the important task of determining the initial number,
viewed as unknown but fixed,  of molecules in a PCR amplification,
i.e. classical quantitative PCR. In literature, it has been treated
under the simplifying assumption of constant replication probabilities
$p(z)$, cf. \cite{Olofsson}, \cite{Vikalo}. For an experimental
approach based on differentiation see \cite{Swillens} and for a
mathematical paper, focussing however on mutations in an abstract
formulation see \cite{Piau}.  Through the use of digital PCR
\cite{Vogelstein} and barcoding \cite{Best}, \cite{Stahlberg} new
possibilities and techniques have been introduced. We hope to be able
to treat such frameworks. The present work should be suitable for
calibration and interpolation of density values in realtime PCR 
\cite{Kubista} in the usual way. Observed values yield
  model parameter estimates. Thus specified,  the model 
  delivers predictions of missing values.

In our setup, the value of $v$ turns out to be
crucial, the cases $0<v<1$ and $v=1$ yielding quite different
situations. If the starting efficiency $v\in (0,1)$, then individual molecules
replicate randomly and essentially independently during an intitial
phase. By branching process theory their number will therefore, to
begin with, grow like the product of a random factor and the famous exponential
population growth. Randomness is therefore an essential part of the
initial conditions of later phases with more of interaction with the
environment but also more of deterministic structure, due to law of
large numbers effects. It is in this sense, the original starting number has
been hidden by a 'veil of uncertainty'. 

 If, on the other hand, $v=1$, the first observable process size can be inverted to yield the starting number.

This phenomenon is what we investigate, for PCR in the present paper
and for populations in habitats with a finite carrying capacity in a
companion paper \cite{CJK}, cf. also \cite{BHKK}, \cite{Barbour}.
 For somewhat related early examples from epidemic processes and a recent from population genetics, cf. 
 \cite{Kendall56,Whittle55,MA15}.

\section{Mathematical setup}

Denote the number of molecules in the $n$-th PCR cycle by $Z_n$,
$n=0,1,2,\ldots$,  so that  $Z_n$ can be viewed as generated by the
recursion
\begin{equation}
\label{Zn}
Z_n =   Z_{n-1} +\sum_{j=1}^{Z_{n-1}} \xi_{n,j},
\end{equation}
started at $Z_0$, where the $\xi_{n,j}$'s are Bernoulli random variables taking values 1 and 0  with complementary probabilities, and
$$
\P\big(\xi_{n,j}=1|Z_{n-1}\big) = \P\big(\xi_{n,j}=1|\F_{n-1}\big) = \frac{vK}{K+Z_{n-1}},
$$
where $\F_{n-1}$ denotes the sigma-algebra of the events, observable before time $n$.

Consider the process $X_n=Z_n/K$, which we shall call the density process. An important role in its  behaviour is played by the function
\begin{equation}\label{fnf}
f(x)=x+\frac{vx}{1+x},
\end{equation}
which is, indeed, the conditional expectation of $X_n$ given $X_{n-1}=x$,
$$
\E(X_n |X_{n-1}=x)=f(x).
$$
The following result is known, see \cite{Kurtz}, \cite{Klethr}.
\begin{theorem} Suppose that $X_0\to x_0$, as $K\to\infty$. Then, for any
  $n$, $$ X_n \xrightarrow[K\to\infty]{\P} f_n(x_0)$$   where $f_n$ denotes the $n$-th iterate of $f$.
 \end{theorem}
If the PCR starts from a fixed number $Z_0$ of molecules, clearly $Z_0/K\to 0$. Since  $f(0)=0$, also  $f_n(0)=0$, for any $n$, and it follows that $\lim_{K\to\infty}X_n= 0,$ for any $n$. In other words, the limiting reaction is not observable at any fixed number of repetitions.
The main result of this paper is that it becomes observable when the number of iterations is $n=\log_b K$, where $b=1+v$.

To arrive at the result we make use of a  linear
replication  process $Y_n$, in which the probability of successful
molecular replication is constant and equal to $v$. In each round each
molecule is thus replaced by two with probability $v$, but remains there
alone with probability $1-v$. The expected number of successors is
thus $1-v+2v=1+v=b$. Mathematically, this process is given recursively by
\begin{equation}
\label{Yn}
Y_n =  Y_{n-1} + \sum_{j=1}^{Y_{n-1}} \eta_{n,j},
\end{equation}
where the $\eta_{n,j}$ are independent Bernoulli random variables with
$$
\P(\eta_{n,j}=1)=v.
$$
Since the $ Y_n/ b^{n}$  constitute  a uniformly integrable
martingale, it has an a.s. limit
\begin{equation}\label{branchlim}
W := \lim_{n\to\infty} b^{-n} Y_n
\end{equation}
with $\E[W] = 1$, provided $Y_0=Z_0=1$.

If the process  starts from  $Z_0$ molecules, then in view of the branching property, the
corresponding limit is $$W(Z_0)=\sum_{i=1}^{Z_0}W_i,$$ where the $W_i$ are i.i.d. with the same continuous 
distribution as $W$.
As is well known from branching process theory (see e.g. Theorem 8.2 in \cite{Harris}), 
the moment generating function of the latter $\phi(s) = \E[e^{-sW}]$, is unique among moment generating 
functions satisfying the functional equation 
$$
\phi(ms)= h(\phi(s)), \quad s\ge 0 
$$
subject to $\phi'(0)=-1$, where $h(s)=\E (s^{Y_1}|Y_0=1)$ and $m=\E (Y_1|Y_0=1)$. In our case, it takes the form
\[\phi((1+v)s) = (1-v)\phi(s) + v\phi(s)^2.\]

The random variable $W(Z_0)$ appears in the main result as an argument of the deterministic function $H$ obtained as
the limit
\begin{equation}\label{fnH}
H(x)=\lim_{n\to\infty}f_n(x/b^n).
\end{equation}
Its existence and some properties are studied in the next section. Here we formulate the main result and an important corollary.

\begin{theorem}\label{main}

Let $v\in (0,1]$ and start the  PCR amplification from $Z_0$ molecules.   Then $X_{\log_{b}K}$ converges in distribution   
$$
X_{\log_{b}K} \xrightarrow[K\to\infty]{ D} H(W(Z_0)),
$$
along any subsequence, such that $\log_{b}K$ are integers.
\end{theorem}
\begin{remark}
With $v=1$, the process $Z_n$ grows deterministically at the geometric rate $b=2$ and in this case $W(Z_0)=Z_0$.
As will be increasingly clear, there are, however reasons to treat $v=1$ separately.
\end{remark}

\begin{corollary}\label{cor1}
For $v\in (0,1]$ and any fixed $n$
\begin{equation}\label{claimbn}
X_{\log_{b}K+n} \xrightarrow[K\to\infty]{ D} f_n(\tilde X_0),
\end{equation}
where   $f_n$ denotes the $n$-th iterate of $f$ and
$$
\tilde X_0=H(W(Z_0)).
$$
This assertion  extends to weak convergence of the sequences  regarded as random elements in 
$\mathbb{R}^{\mathbb{Z}}$:
$$
\{X_{\log_{b}K+n}\}_{-\infty}^{\infty} \xrightarrow[K\to\infty]{ D} \{f_n(\tilde X_0)\}_{-\infty}^{\infty}.
$$

\end{corollary}

\begin{remark}
The limits increase strictly with respect to $n$. If $0<v<1$, their entries are continuous random variables with positive variance, 
whereas if $v=1$ they are positive reals. If the limit in \eqref{claimbn} is taken along an arbitrary subsequence $K$, then
$
X_{[\log_{b}K]}
$
is asymptotic to the same limit up to a deterministic correction, which emerges in the rounding:
$$
X_{[\log_{b}K]} - H\Big(W(Z_0)b^{[\log_b K]-\log_b K}\Big)\xrightarrow[K\to\infty]{D} 0.
$$
\end{remark}

\section{The limit function $H(x)$}

\subsection{Existence}
Write the two expressions for $f$, \eqref{fnf} and
\begin{equation}\label{near0}
f(x)=bx-\frac{v x^2}{1+x}=bx-g(x),
\end{equation}
 where   $g(x)=\frac{v x^2}{1+x}$. This expression is more suitable for analysis of iterates of $f$ near zero.

 It is easy to establish that $f$ is increasing,  which yields that all $f_n$ are increasing. Since $g(x)> 0$ for any $x>0$,
 $$f(x/b)<x.$$
Hence
$$
f_{n+1}(x/b^{n+1})=f_n(f(x/b^{n+1}))<f_n( x/b^{n}),
$$
and the sequence $f_n(x/b^n)$ is monotone decreasing in $n$ for any positive $x$.
Therefore the following limit in \eqref{fnH} exists,
$$
H(x)=\lim_{n\to\infty}f_n(x/b^n).
$$

\subsection{Continuity} 
We show next that the convergence in \eqref{fnH} is uniform on bounded intervals. First observe that
$$f'(x)=1+\frac{v}{(1+x)^2}\le 1+v=b.$$
It is now easy to see by induction, that for any $n$ and $x$
$$f_n'(x)\le b^n.$$
Next, by \eqref{near0} the Taylor expansion reads
$$
f_{n+1}(x)=f_n(f(x))=f_n(bx-g(x))=f_n(bx)-f_n'(\theta_n)g(x),
$$
for an appropriate $\theta_n$. Replace now $x$ by $x/b^{n+1}$ to have
$$f_{n+1}(x/b^{n+1})= f_n(x/b^n)-f_n'(\theta_n)g(x/b^n).$$
Hence we obtain
\begin{equation}\label{fund1}
f_n(x/b^n)-f_{n+1}(x/b^{n+1})=f_n'(\theta_n)g(x/b^n)\le b^n g(x/b^n)\le vx^2b^{-n},
\end{equation}
where we have used that $g(x)=vx^2/(1+x)\le vx^2$.
The bound \eqref{fund1} shows that the series $$\sum_{n=0}^\infty f_{n+1}(x/b^{n+1})-f_n(x/b^n)$$ converges uniformly 
on compacts. As a consequence of uniform convergence, we have that $H$ is continuous.

\subsection{The functional equation}

Further, since $f_{n+1}(x/b^{n+1})=f(f_n((x/b)/b^{n}))$,  by taking the limit as $n\to\infty$, we obtain that $H$ solves 
Schr\"{o}der's functional equation
\begin{equation}\label{AbelH}
H(x)=f(H(x/b)).
\end{equation}
However, since the zero function is a solution, we must show that $H$ is not identically zero. $H(x)=\infty$ is also a solution,  it is however directly excluded, since convergence is from above, $f_n(x/b^n)>H(x)$.

To show that $H$ is positive, use \eqref{near0} to obtain
the following formula for the $n$-th iterate
$$
f_n(x)=b^nx-\sum_{i=0}^{n-1} b^{n-1-i}g(f_{i}(x)),
$$
where, as usual, $f_0(x)=x$.
Replacing $x$ with $xb^{-n}$, we have
\begin{equation}\label{fnbn}
f_n(xb^{-n})= x-\sum_{i=0}^{n-1} b^{n-1-i}g(f_{i}(xb^{-n})).
\end{equation}
Clearly, $f_i(x)\le b^i x$, and $g(x)\le v x^2$, therefore
$$
b^{n-1-i}g(f_{i}(xb^{-n}))\le v  b^{n-1-i} (b^i xb^{-n})^2=v x^2 b^{-n+i-1},
$$
 and
$$
\sum_{i=0}^{n-1} b^{n-1-i}g(f_{i}(xb^{-n}))\le v x^2  \sum_{i=0}^{n-1} b^{-n+i-1}\le x^2.
$$
Hence from \eqref{fnbn}, for any $n$
$$f_n(xb^{-n})\ge x-x^2,$$
which is strictly positive for $0<x<1$. Therefore  $H(x)>0$ in this domain.

\subsection{Monotonicity}

We show next that $H$ is increasing. Let $H_n(x)=f_n(x/b^n)$. Then each $H_n(x)$ is increasing
and thus $H(x)=\lim_{n\to\infty}H_n(x)$ does not decrease. Further, recall that
$$
f'(x)=1+\frac{v}{(1+x)^2} = b-vx \frac{ 2 +x }{(1+x)^2} > b-2x  
$$
and $f_j(x/b^j) \le  x$ for all $j\ge 0$.  Hence for any $x\le b^2/2$, 
\begin{align*}
H_n'(x) & =b^{-n} f'_n(x/b^n) = b^{-n}\prod_{j=0}^{n-1} f'(f_j(x/b^n)) \ge 
b^{-n}\prod_{j=0}^{n-1} \big(b- 2f_j(x/b^n)\big)\\
&\ge 
b^{-n}\prod_{j=0}^{n-1} \big(b- 2x b^{j-n}\big)
\ge \prod_{j=0}^{n-1} \big(1-  b^{-j}\big)\ge e^{-v}, \quad \forall \, n\ge 0,
\end{align*}
and 
$$
H_n(x_2)-H_n(x_1) = \int_{x_1}^{x_2} H'_n(x)dx > (x_2-x_1) e^{-v}>0, \quad x_1<x_2<b^2/2.
$$
Taking the limit $n\to\infty$, we see that $H(x)$ is a strictly increasing function on an open vicinity of the origin. 

Suppose now that $H$ is constant on an interval $[x_1,x_2]$ with $x_2>x_1$. Then, by \eqref{AbelH},
$
H(x_1/b^k) = H(x_2/b^k)
$
for any integer $k\ge 1$ and, since $H(x)$ does not decrease, it must be constant on all the intervals $[x_1/b^k,x_2/b^k]$. 
In particular, $H(x)$ cannot be strictly increasing on any open vicinity of the origin. The obtained contradiction shows
that $H$ is strictly increasing everywhere on $\Real_+$.

Next, since we have shown that the $H_n$ converge uniformly,
 $$H_n(x+o_n(1))\to H(x),$$ for any $o_n(1)\to 0$ as $n\to\infty$. Thus we have the following corollary needed
in the proofs to come.
\begin{corollary}\label{olittle}
$$\lim_{n\to\infty}f_n(x/b^n+o(b^{-n}))=H(x).$$
\end{corollary}

We shall also need the inverse $G:=H^{-1}$.
It is easy to see that it solves the functional equation
$$
G(x)=\frac{1}{b}G(f(x)).
$$

%

 \section{Proofs}


Let us start with the fundamental recursive equation for the stochastic density process $X_n$ (cf. \cite{Klethr})
\begin{equation}
\label{barZ}
X_n = f(X_{n-1})  + \frac 1{\sqrt K}\eps_{n},
\end{equation}
with
$$
\eps_n =  \frac 1{\sqrt K}\sum_{j=1}^{KX_{n-1}} (\xi_{n,j}-E(\xi_{n,j}|\F_{n-1})).
$$
Note that $\eps_n$ is a martingale difference sequence
 $\E(\eps_n|\F_{n-1})=0$  and
\begin{equation}
\label{varepsn}
\E (\eps^2_n|\F_{n-1}) =   \frac{vX_{n-1} }{1+X_{n-1}}\left(1- \frac{ v  }{1+X_{n-1}}\right)\le v.
\end{equation}

The corresponding deterministic recursion, obtained by omitting the martingale difference term, is
\begin{equation}\label{xn}
x_n=f(x_{n-1})=f_n(x_0).
\end{equation}

\subsection{Proof  of Theorem \ref{main}} In what follows bar denotes the density processes, i.e., 
$\bar Z_n = Z_n/K$, $\bar Y_n = Y_n/K$. Consider first the case $v<1$.
Define times
$$
n_1 = c \log_b K\quad \text{and} \quad n_2 = \log_b K,
$$
where $c\in (\frac 1 2,1)$ is an arbitrary fixed constant and $K$ is such that both $n_1$ and $n_2$ are integers.

The crux of the proof is to approximate the density process $X_n=\bar Z_n := Z_n/K$ in two steps.
First, on the interval $[0,n_1]$ by the linear process $\bar Y$, and then  on the interval $[n_1,n_2]$ by the nonlinear deterministic recursion, however started from the random point $\bar Y_{n_1}$, resulting from the first step.

Denote by $\phi_{k,\ell}(x)$ the flow, generated by the nonlinear deterministic recursion \eqref{xn}, i.e. its solution at time 
$\ell$, when started from $x$ at time $k$,  $x_\ell=\phi_{k,\ell}(x_k)=f_{l-k}(x_k)$.
Further, write $\Phi_{k,\ell}(x)$ for the stochastic flow generated by the nonlinear process $X$, that is, 
the random map defined by the solution of the equation, cf. \eqref{Zn}, 
$$
X_n =   X_{n-1} +\sum_{j=1}^{K X_{n-1}} \xi_{n,j}, \quad n = k+1,...,\ell
$$
subject to $X_k=x$, at the terminal time $n:=\ell$.  In particular, $X_k = \Phi_{k,\ell}(X_\ell)$ for any $k> \ell > 0$,
and
\begin{multline*}
X_{n_2}    =  \Phi_{n_1,n_2}(X_{n_1})=\phi_{n_1,n_2}(X_{n_1})+ (\Phi_{n_1,n_2}(X_{n_1})-\phi_{n_1,n_2}(X_{n_1})) = \\
   \phi_{n_1,n_2}(\bar Y_{n_1})    +  (\Phi_{n_1,n_2}(X_{n_1})-\phi_{n_1,n_2}(X_{n_1}))
    +  (\phi_{n_1,n_2}(X_{n_1})-\phi_{n_1,n_2}(\bar Y_{n_1})).
\end{multline*}
Let us stress that all the random objects here are defined on the same probability space and by construction coupled 
as described at the beginning of the proof. 
  
In the next steps we show that
\begin{equation}\label{term1}
 \phi_{n_1,n_2}(\bar Y_{n_1})\xrightarrow[K\to\infty]{\text{a.s.}} H(W(Z_0)),
\end{equation}
\begin{equation}
\label{dva}
\Phi_{n_1,n_2}(X_{n_1})-\phi_{n_1,n_2}(X_{n_1})\xrightarrow[K\to\infty]{\P} 0,
\end{equation}
and
\begin{equation}
\label{raz}
\phi_{n_1,n_2}(X_{n_1})-\phi_{n_1,n_2}(\bar Y_{n_1})\xrightarrow[K\to\infty]{\P} 0.
\end{equation}


By  \eqref{branchlim}, with $W=W(Z_0)$, we may write
$$
Y_{n_1}=Wb^{n_1}+o(b^{n_1})=Wb^{c\log_b K}+o(b^{c\log_b K}),
$$
and hence 
$$\bar Y_{n_1}=\frac{1}{K}Y_{n_1}=Wb^{-(1-c)\log_b K}+o(b^{-(1-c)\log_b K}).$$
Therefore,  \eqref{term1} follows from Corollary \ref{olittle}, 
\begin{align*}
&
 \phi_{n_1,n_2}(\bar Y_{n_1})  =f_{n_2-n_1}(\bar Y_{n_1}) \\
 & = f_{(1-c)\log_bK}(Wb^{-(1-c)\log_b K}+o(b^{-(1-c)\log_b K}))\xrightarrow[K\to\infty]{\text{a.s.}} H(W).
\end{align*}
To show \eqref{dva}   let for $n>n_1$ $$\delta_n = \E |\Phi_{n_1,n }(X_{n_1})-\phi_{n_1,n }(X_{n_1})|.$$

Subtracting the deterministic recursion \eqref{xn} from the stochastic one \eqref{barZ}    we have
$$
X_n - x_n =X_{n-1}-x_{n-1}+
v\frac{  X_{n-1}-x_{n-1} }{\big(1+X_{n-1}\big)\big(1+x_{n-1}\big)} + \frac 1{\sqrt K}\eps_{n}.
$$

Thus the sequence $\delta_n$  satisfies
$$
 \delta_n  \le
b  \delta_{n-1}    + \frac 1{\sqrt K} \sqrt v,
$$
 where we have used \eqref{varepsn} to bound $\E |\eps_{n}|$. Note that   $\delta_{n_1}=0$, as both recursions start  at the same point $X_{n_1}$ at time $n_1$.
Therefore
$$
 \delta_{n_2}  \le \sqrt v\frac 1{\sqrt K}\sum_{j=0}^{n_2-n_1-1}b^{j}\le C  K^{-\frac 1 2} b^{n_2-n_1}\le  C     K^{\frac 1 2-c}\xrightarrow[K\to\infty]{}0,
$$
since $c>\frac 1 2$ and \eqref{dva} now follows.

The proof of  \eqref{raz} is more delicate and  is done by coupling.
We  construct the nonlinear and linear replication processes $Z_n$ and $Y_n$ on the same probability space as follows. Let $U_{n,j}$ $n,j\in \mathbb{N}$
be i.i.d. random variables with the uniform distribution on $[0,1]$. Define
$$
\xi_{n,j} = \oneset{\left\{U_{n,j}\le \frac{vK}{K+Z_{n-1}}\right\}} \quad \text{and} \quad \eta_{n,j} = \one{U_{n,j}\le v}.
$$
Then $Z_n$ and $Y_n$ are realized by the formulae \eqref{Zn} and \eqref{Yn} with $\xi_{n,j}$ and $\eta_{n,j}$ as above. 
Since $\frac{vK}{K+Z_{n-1}}<v$, we have $\xi_{n,j}\le \eta_{n,j}$ for all $n,j$ and  therefore the linear process $Y$ is always greater than the nonlinear process $Z$, $$Z_n\le Y_n,\;\mbox{ for all}\; n.$$ Construct an auxilliary linear process $V_n$, which bounds $Z_n$ from below until $Z_n$ gets larger than
$K^\gamma$ for  $\gamma\in (0,1)$. Actually we require that $c<\gamma <1$. Let
$$
\zeta_{n,j} = \oneset{\left\{U_{n,j}\le \frac{vK}{K+K^\gamma}\right\}},
$$
and
$$V_n=V_{n-1}+\sum_{j=1}^{V_{n-1}}\zeta_{n,j}.$$
Then clearly,  $\zeta_{n,j}<\xi_{n,j}$ as long as $Z_{n-1}<K^\gamma$.  Hence
$$V_n\le Z_n,\;\mbox{for}\;  n< \tau=\inf\{k: Z_k> K^\gamma\}.$$
It is also clear that for all $n,j$, $\zeta_{n,j}<\eta_{n,j}$
hence $V_n\le Y_n$. Thus we obtain
\begin{equation}\label{YminusZn1}
\begin{aligned}
Y_n-Z_n&= Y_n-V_n+V_n-Z_n \\
&\le  Y_n-V_n+(V_n-Z_n)1_{n> \tau} \\
&\le  Y_n-V_n+V_n1_{\tau<n}.
\end{aligned}
\end{equation}
We show next that
\begin{equation}\label{Kclimit}
\lim_{K\to\infty}(Y_{n_1}-Z_{n_1})K^{-c} =0
\end{equation}
by using the inequality above. Since the moments of simple Galton-Watson processes are easily computed (Theorem 5.1 in \cite{Harris})
$$
\E V_{n_1}=(1+\frac{v}{1+K^{\gamma-1}})^{c\log_bK}=b^{c\log_bK}(1-\frac{v}{b(1+K^{\gamma-1})}K^{\gamma-1})^{c\log_bK}\sim K^c.
$$
Since $\E Y_{n_1}=b^{n_1}=K^c$ also, the first term in \eqref{YminusZn1} satisfies
$$\lim_{K\to\infty}\E (Y_{n_1}-V_{n_1}) K^{-c}=0.$$
By the Cauchy-Schwartz inequality for the second term
$$\E V_{n_1}1_{\tau<n_1}\le \Big(\E V^2_{n_1}\P (\tau<n_1)\Big)^{1/2}.$$
Since $Z_n<Y_n$ for all $n$, it takes longer for the former process to reach $K^\gamma$  than  the corresponding time for the latter,
$$\tau\ge \sigma=\inf\{n:Y_n>K^\gamma\}.$$
Therefore
\begin{align*}
\P (\tau<n_1)&\le  \P \big(\sigma <n_1\big)\\
&=  \P (\sup_{n <n_1}Y_n>K^\gamma)\le \P (b^{-n_1}\sup_{n <n_1}Y_n>K^\gamma b^{-n_1})\\
&\le  \P (\sup_{n <n_1}Y_nb^{-n}>K^{\gamma -c})\le K^{c-\gamma},
\end{align*}
where the last bound is Doob's inequality for the martingale $Y_nb^{-n}$. Taking into account that
$\E V^2_{n_1}\sim K^{2c}$, we obtain from the above estimates
$$
\lim_{K\to\infty} K^{-c}\E V_{n_1}1_{\tau<n_1}=0.
$$
Recall that $\gamma > c$. It follows that the convergence to the limit in \eqref{Kclimit} holds in $L^1$, and in probability.
For the corresponding densities, we have by dividing through by $K$ that
\begin{equation}\label{Kclimdens}
\lim_{K\to\infty}(\bar Y_{n_1}-X_{n_1})K^{1-c} =0
\end{equation}
Since $\phi_{n_1,n_2}(x)=f_{n_2-n_1}(x)$ and the function $f$ is concave ($f''<0$), its derivative attains its maximum vaue at zero,  $f'(0)=b$ and $f'_n(x)\le b^n$ for any $x\geq 0$.
Therefore   $|f_n(x)-f_n(y)|\le b^n |x-y|$. For $y=\bar Y_{n_1}$ and $x=X_{n_1}$, this and \eqref{Kclimdens} yields
\begin{eqnarray*}
0\le f_{n_2-n_1}(\bar Y_{n_1})-f_{n_2-n_1}(X_{n_1})&\le& b^{n_2-n_1}(\bar Y_{n_1}-X_{n_1})\\
&=&K^{1-c}(\bar Y_{n_1}-X_{n_1})\to 0,
\end{eqnarray*}
and the proof of case $v<1$ is complete.

Consider now the case $v=1$.
In this case, the probability of successful replication is
$$
\P\big(\xi_{n,j}=1|Z_{n-1}\big) = \frac{K}{K+Z_{n-1}},
$$
and the function $f$ is
$$
f(x)=x+\frac{x}{1+x}.
$$
Here $b=v+1=2$ and
$$
H(x)=\lim_{n\to\infty}f_n(x/2^n).
$$
The proof is the same, except that   the linear replication process $Y_n$ is in fact deterministic $Y_n=Z_02^n$,  if it starts with $Z_0$ molecules, because the probability of replication is 1, $\P (\eta_{n,j}=1)=v=1$.
Hence the limit $W=Y_n/2^n=Z_0$. The theorem is proved.

\subsection{Proof  of Corollary \ref{cor1}}

The result follows by induction on $n$ from the fundamental representation \eqref{barZ}. For $n=0$ it is the statement of the main result. For $n=1$ take limits as $K\to \infty$ in \eqref{barZ}, and note that the stochastic term vanishes. Similarly, having proved it for $n$, it follows for $n+1$. The functional limit theorem follows from 
finite dimensional convergence implying convergence in the sequence space, {\em cf.} \cite{B}, p. 19.

\section{The relation to actual observations}

Let $\rho$ denote the minimal observable concentration of DNA in the PCR experiment under consideration. Assume that the latter starts from $z=Z_0$ inititial templates, where $z$ is an unknown number and $x=X_0=z/K < \rho$. Our aim is to determine $z$ for $K >> z$. Mathematically, we shall interpret this as $K\to\infty$. In PCR literature based on enzyme kinetic considerations, values of the Michaelis-Menten constant range at least from $10^6$ \cite{lalam06} up to $10^{15}$ \cite{Gevertz}, in terms of molecule numbers.

There are then two cases, known or unknown rate $v$. In the latter situation, $v$ will have to be estimated from the observed concentrations. Further, as pointed out, the cases $v=1$ and $v<1$ exhibit an intriguing disparity, viz. consider first $v<1$. By Corollary \ref{cor1}
$$
\big\{X_{\log_{b}K+n}\big\}_{-\infty}^{\infty} \xrightarrow[K\to\infty]{ D} \big\{f_n(H(W(z)))\big\}_{-\infty}^{\infty}.
$$
The limit process here has strictly increasing trajectories and its entries have continuous distributions, so with probability one 
none of them equals $\rho$. The first hitting time 
$$
(x_n)\mapsto \inf\big\{n\in \mathbb{Z}; x_n\geq \rho\big\}, \quad x\in \Real^{\mathbb{Z}}
$$
being a discontinuous functional with respect to the locally uniform metric on space of sequences, 
is however continuous almost surely under the limit law. Therefore 
$$
\tau^K(\rho) := \inf\big\{n\in\mathbb{Z};X_{\log_{b}K+n}\geq \rho\big\}
$$ 
converges weakly 
to 
$$
\tau(\rho) := \inf\{n\in\mathbb{Z};f_n(H(W(z)))\geq \rho\} \quad \text{as }  K\to\infty.
$$

If $v=1$, the limit sequence is deterministic and strictly increasing. Provided no $f_n(H(z))$ happens to coincide with $\rho$, we have weak convergence $\tau^K(\rho)\to \tau(\rho)$. Otherwise, $\lim_{K\to\infty}\tau^K(\rho)$ still exists and differs at most by 1 from $\tau(\rho)$. 

\begin{figure*}
\input{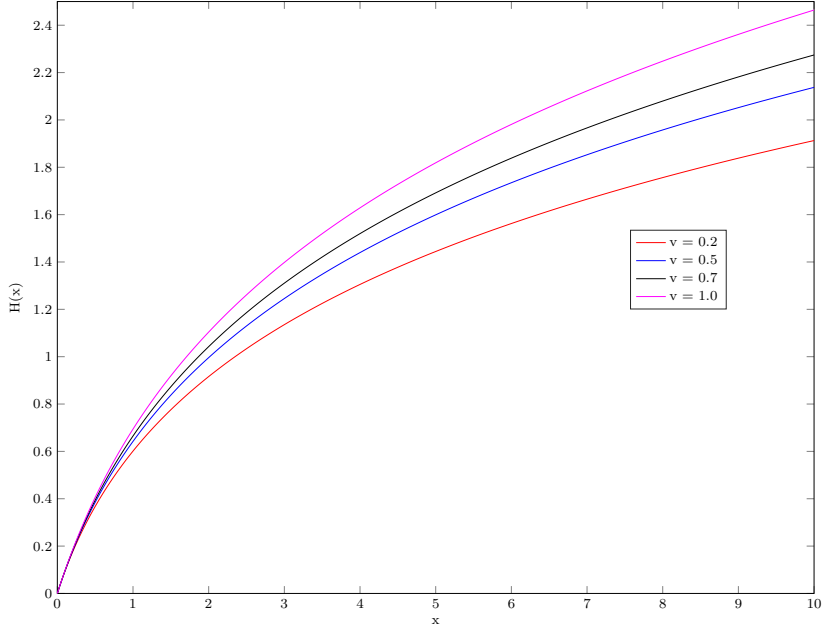} 
\caption{The function $H(x)$ for several values of $v$}
\label{fig1} 
\end{figure*}

We disregard this  nuisance and assume in both cases that we have observed concentration values strictly larger than $\rho$ from $ \log_{b}K+\tau^K(\rho)\approx  \log_{b}K + \tau(\rho)$ onwards: $\kappa_0=f_\tau(H(W(z))),\kappa_1=f_{\tau+1}(H(W(z)),\kappa_2= f_{\tau+2}(H(W(z)), \ldots$, and correspondingly for $v=1$,  
$\kappa_0=f_\tau(H(z))$, $\kappa_1=f_{\tau+1}(H(z))$, $\kappa_2= f_{\tau+2}(H(z)), \ldots$ (to ease notation, we omit the dependence of $\tau$ upon $\rho$.)  By \eqref{AbelH} this simplifies to
\[\kappa_j= H(W(z)b^{\tau+j})\]
for $v<1$ and
\[\kappa_j= H(zb^{\tau+j})\]
otherwise. Note that typically, since the experimenter would like to catch the density as early as possible,  $\kappa_0 \approx \rho$, which for example could be of the order of 0.05. Since $H(x)$ is fairly close to the diagonal $H(x)=x$ for $0\leq x\leq 0.5$ (see Figure \ref{fig1}) and $W(z)\approx z$, we can conclude that as a rule $\tau<0$. 

As well as assuming $K$ and $\rho$ known it is easy to think of situations where so is $v$. Then we can proceed directly to determining $z$. For $v=1$ this is straightforward:
$$z = b^{-\tau}G(\kappa_0).$$
More generally,
$$z = b^{-\tau-j}G(\kappa_j).$$
If there is variation between the $z$-values thus obtained we can of course take arithmetic means of the right hand side for the different observed $j$.

Now, if $v<1$, we obtain
$$ \sum_{i=1}^zW_i= W(z) =  b^{-\tau}G(\kappa_0),$$
in the sense that the right hand side is an observed value of the random variable $W(z)$. The initial number $z$ of DNA molecules has now been hidden from direct calculation. What can be done is to estimate $z$ from data, e.g. maximise the density at the first point of observation,
$$\psi^{*z}(b^{-\tau}G(\kappa_0)),$$
where * denotes convolution power, $\psi$ is the density of $W$, which we know to have the moment generating function $\phi$ from Section 2, corresponding to $v$. In this, $z$ is an unknown parameter and we obtain  a  maximum likelihood estimate 
$\hat{z}= \mathrm{argmax}_z\psi^{*z}(t)$, where $t= b^{-\tau}G(\kappa_0)$ and $z$ ranges over natural numbers. Again we can also consider later $\kappa$-values and take averages, if this increases stability. Note that if $z$ is large (but still much smaller than $K$), then by the local central limit theorem the ML problem is roughly the same as finding $z$ maximizing the normal density with mean $z$ 
and variance $z\frac{1-v}{1+v}=:z\sigma^2$ 
at the point  $t= b^{-\tau}G(\kappa_0)$,
$$\phi^{*z}(t)\approx \sqrt{\frac{1+v}{2\pi z(1-v)}}\exp\frac{-(t-z)^2}{2z(1-v)/(1+v)}.$$
This yields the estimate 
$$
\hat{z} =\sqrt{t^2 + \sigma^4/4} - \sigma^2/2 =
\sqrt{\big(b^{-\tau}G(\kappa_0)\big)^2 - \frac 1 4\left(\frac{1-v}{ 1+v }\right)^2} - \frac 1 2 \left(\frac{1-v}{ 1+v }\right)^2,
$$ 
or rather one of its neighboring integers. 

Now, if entities cannot be deduced {\em a priori} the question arises to what extent they can be estimated from our sequence of observations. Clearly, in the limit the relation between an observation $x$ and its successor in the next round will be that the latter converges to $f(x)$, as $K\to\infty$, by Corollary \ref{cor1}. Thus e.g.,
\[\kappa_1 =\kappa_0 + \frac{v\kappa_0}{1+\kappa_0}\]
or 
\[v= \kappa_1(1+\kappa_0) - 1.\]

These problems are fairly standard in statistical literature
  but certainly deserve a special investigation in the present
  context, if possible together with an experimental study of
  replication of single or few molecules, in order to determine the
  initial efficiency, $v$.


\end{document}